\documentclass[preprint,1p]{elsarticle}

\usepackage{amssymb}
\usepackage{amsthm}
\usepackage{amsmath}
\usepackage{geometry,graphicx,amsfonts,bm,url}

\usepackage{natbib}

\usepackage{latexsym,bm,amssymb,amsfonts,amsbsy,amsmath,graphics,graphicx,color}
\usepackage{pifont}
\usepackage{url}
\usepackage{amsmath}
\usepackage{amsfonts}
\usepackage{amssymb}

\newtheorem{theorem}{Theorem} 
\newtheorem{corollary}[theorem]{Corollary}

\def\bd{\begin{description}}
	\def\ed{\end{description}}

\def\beq{\begin{equation}}
\def\eeq{\end{equation}}
\def\bea{\begin{eqnarray}}
\def\eea{\end{eqnarray}}
\def\beas{\begin{eqnarray*}}
\def\eeas{\end{eqnarray*}}

\def\RR{\mathbb R}

\def\pmatrix{\left(\begin{array}}
\def\endpmatrix{\end{array}\right)}

\def\one{{\bf 1}}

\newlength\figureheight
\newlength\figurewidth

\newcommand{\be}{\begin{equation}}
\newcommand{\ee}{\end{equation}}

\usepackage{graphicx,color,url}

\usepackage{lipsum}
\usepackage{amsfonts}
\usepackage{graphicx}
\usepackage{epstopdf}
\usepackage{algorithmic}
\ifpdf
\DeclareGraphicsExtensions{.png,.eps,.pdf,.jpg}
\else
\DeclareGraphicsExtensions{.eps}
\fi

\usepackage{amsopn}

\journal{Applied Mathematics Letters}
\begin{document}
	
\begin{frontmatter}

\title{High-order Gauss-Legendre methods admit a composition representation and a conjugate-symplectic counterpart}

\author[label0]{Felice Iavernaro}
\affiliation[label0]{organization={Dipartimento di Matematica, Università  degli Studi di Bari Aldo Moro},
	addressline={Via Orabona 4}, 
	city={Bari},
	postcode={70125}, 
	country={Italy}}

\author[label1]{Francesca Mazzia} 
\affiliation[label1]{organization={Dipartimento di Informatica,  Università degli Studi di Bari Aldo Moro},
	addressline={Via Orabona 4}, 
	city={Bari},
	postcode={70125}, 
	country={Italy}}

\author[label2]{Ernst Hairer}
\affiliation[label2]{organization={Section de Mathématiques, Université de Genève}, 
	addressline={24 rue du Général-Dufour}, 
	city={Genève},
	postcode={1211 Genève 4}, 
	country={Switzerland}}

\begin{abstract}
One of the most classical pairs of symplectic and conjugate-symplectic schemes is given by the Midpoint method (the Gauss--Runge--Kutta method of order 2) and the Trapezoidal rule. These can be interpreted as compositions of the Implicit and Explicit Euler methods, taken in direct and reverse order, respectively. This naturally raises the question of whether a similar composition structure exists for higher-order Gauss--Legendre methods. In this paper, we provide a positive answer by first examining the fourth-order case and then outlining a generalization to higher orders. 
\end{abstract}

\end{frontmatter}

\section{Introduction and main result}
\label{sec:intro}
In the context of the numerical solution of initial value problems, it is well known that composing the Explicit and Implicit Euler methods (denoted EE and IE, respectively) with halved stepsize gives rise to two fundamental numerical integrators: the \emph{trapezoidal rule} (TR), obtained by composing EE followed by IE, and the \emph{implicit midpoint rule} (MP), obtained by reversing their order:
$$
\text{TR} = \text{IE} \circ  \text{EE} , \qquad \text{MP} = \text{EE} \circ  \text{IE} .
$$
Both methods are of order two and exhibit remarkable structure-preserving properties when applied to Hamiltonian systems of the form
\begin{equation}\label{ham}
	y'(t) = f(t,y) := J \nabla H(y(t)), \qquad y(t_0) = y_0 \in \RR^{2m}, \qquad
	y = \pmatrix{c} q \\ p \endpmatrix, \quad q, p \in \RR^m,
\end{equation}
where $H: \RR^{2m} \to \RR$ is the Hamiltonian function, and $J$ is the canonical symplectic matrix
$$
J = \pmatrix{rr} 0 & I \\ -I & 0 \endpmatrix,
$$
with $I$ denoting the $m \times m$ identity matrix. In this setting, the implicit midpoint method is symplectic, whereas the trapezoidal rule is its conjugate-symplectic counterpart.

Symplectic  integrators enjoy important geometric properties, such as volume preservation of closed surfaces in phase space, conservation of all quadratic first integrals, and near-conservation of the Hamiltonian function over exponentially long time intervals \cite{HLW06,BI16}.  By their very nature, conjugate-symplectic integrators inherit the same qualitative behavior \cite{HZ12}.

A prominent family of symplectic implicit Runge--Kutta methods is given by the Gauss--Legendre collocation schemes, whose second-order member coincides with the implicit midpoint rule. This naturally raises the question: \emph{does a similar composition structure exist for higher-order Gauss--Legendre Runge--Kutta methods?} In particular, one may ask whether these methods can also be represented as compositions of suitable high-order extensions of the Euler schemes, and whether such compositions admit meaningful conjugate-symplectic counterparts.

In this paper, we first address this question for the two-stage, fourth-order Gauss--Legendre method, defined by the tableau 
$$
\begin{array}{c|cc}
	\frac{1}{2} - \frac{\sqrt{3}}{6}  & \frac{1}{4} & \frac{1}{4} - \frac{\sqrt{3}}{6} \\[.15cm]
	\frac{1}{2} + \frac{\sqrt{3}}{6}  & \frac{1}{4} + \frac{\sqrt{3}}{6} & \frac{1}{4} \\[.15cm]
	\hline
	& \frac{1}{2}  & \frac{1}{2} 
\end{array}
$$
By exploiting discrete variants of two multi-derivative high-order one-step formulae (see \cite{Iavernaro2018_torino, Mathematica2021}), we construct two two-stage Runge--Kutta methods, denoted $\Phi_h$ and $\Psi_h$, with the following Butcher tableau:
\begin{equation}
	\label{PhiPsi}
	\Phi_h := \quad
	\begin{array}{c|cc}
		1 - \frac{\sqrt{3}}{3}  & \frac{1}{2} & \frac{1}{2} - \frac{\sqrt{3}}{3} \\[.15cm]
		1 + \frac{\sqrt{3}}{3}  & \frac{1}{2} + \frac{\sqrt{3}}{3} & \frac{1}{2} \\[.15cm]
		\hline
		& \frac{1}{2} + \frac{\sqrt{3}}{4} & \frac{1}{2} - \frac{\sqrt{3}}{4}
	\end{array}
	\hspace{1.5cm}
	\Psi_h := \quad
	\begin{array}{c|cc}
		- \frac{\sqrt{3}}{3} & -\frac{\sqrt{3}}{4} & -\frac{\sqrt{3}}{12} \\[.15cm]
		\frac{\sqrt{3}}{3}  & \frac{\sqrt{3}}{12} & \frac{\sqrt{3}}{4} \\[.15cm]
		\hline
		& \frac{1}{2} - \frac{\sqrt{3}}{4} & \frac{1}{2} + \frac{\sqrt{3}}{4}
	\end{array}
\end{equation}
It turns out that the composition $\Psi_{h/2} \circ \Phi_{h/2}$ yields the Gauss--Legendre method of order four, while the reverse composition $\Phi_{h/2} \circ \Psi_{h/2}$ gives rise to a novel conjugate-symplectic Runge--Kutta method with the following  tableau:
\begin{equation}\label{PhiPsiConjSymp}
	\begin{array}{c|cccc}
		- \frac{\sqrt{3}}{6}  & - \frac{\sqrt{3}}{8} & - \frac{\sqrt{3}}{24}  &   0 &   0 \\[.15cm]
		+ \frac{\sqrt{3}}{6}     &  \frac{\sqrt{3}}{24} &  \frac{\sqrt{3}}{8}   & 0 &  0\\[.15cm]
		1 - \frac{\sqrt{3}}{6}  & \frac{1}{4}- \frac{\sqrt{3}}{8} & \frac{1}{4} +\frac{\sqrt{3}}{8}  & \frac{1}{4} & \frac{1}{4} - \frac{\sqrt{3}}{6} \\[.15cm]
		1 + \frac{\sqrt{3}}{6}  & \frac{1}{4}- \frac{\sqrt{3}}{8} & \frac{1}{4} +\frac{\sqrt{3}}{8}  & \frac{1}{4} + \frac{\sqrt{3}}{6} & \frac{1}{4} \\[.15cm]
		
		\hline
		& \frac{1}{4} - \frac{\sqrt{3}}{8} & \frac{1}{4} + \frac{\sqrt{3}}{8} & \frac{1}{4} + \frac{\sqrt{3}}{8} & \frac{1}{4} - \frac{\sqrt{3}}{8}
	\end{array}
\end{equation}

	Since for the Dahlquist test equation the application of Runge--Kutta methods commutes, both methods 
	(the Gauss--Legendre method of order 4 and the method  (\ref{PhiPsiConjSymp}) have the same stability function.
Moreover, although both $\Phi_h$ and $\Psi_h$ are implicit, their repeated composition requires solving only one nonlinear system at each time step. 

By analyzing the relationship between the coefficients defining the methods $\Phi_h$ and $\Psi_h$ defined in \eqref{PhiPsi}, we were able to extend the construction to Gauss–Legendre methods of arbitrarily high order.

The paper is organized as follows. In Section~\ref{sec: MDMP_MDTR}, we introduce the multiderivative methods and their discrete variants obtained by replacing derivatives with suitable finite-difference formulae. In Section~\ref{MDMP}, we show that, for an appropriate choice of the discretization parameter, the Gauss method is recovered. Section~\ref{MDTR} introduces the corresponding conjugate scheme.  In Section~\ref{sec:gen}, we present the generalization to higher-order integrators, while Section~\ref{sec:conc} contains our concluding remarks.

\section{Background}
\label{sec: MDMP_MDTR}

Our focus is on two conjugate classes of fourth-order multi-derivative one-step methods: the multi-derivative midpoint (MDMP) and trapezoidal (MDTR) schemes \cite{Iavernaro2018_torino,Mathematica2021}. Denoting by ET4 and IT4 the fourth-order multi-derivative extensions of the explicit and implicit Euler methods, respectively, we define:

\begin{center}
	\begin{tabular}{llll}
		MDMP4 :=~ ET4 $\circ$ IT4 :  \\
		\quad $\displaystyle y_1=y_0+hf(y_{1/2} ) +  \frac{h^3}{24}D_2f(y_{1/2}),$ \\ $\Longleftrightarrow$ 
		
		$ \left\{
		\begin{array}{lll}\displaystyle
			y_{1/2}=y_0+\frac{h}{2}f(y_{1/2}) - \frac{h^2}{8} D_1 f(y_{1/2})  + \frac{h^3}{48}D_2f(y_{1/2}), \\[.25cm]
			\displaystyle y_{1}=y_{1/2}+\frac{h}{2}f(y_{1/2} ) + \frac{h^2}{8} D_1f(y_{1/2})  + \frac{h^3}{48}D_2f(y_{1/2}),  
		\end{array} \right.
		$ 
		\\
		\\
		MDTR4 := ~IT4 $\circ$ ET4:  \\
		\quad $\displaystyle y_1=y_0+\frac{h}{2}\left(f(y_1)+f(y_0)\right) -  \frac{h^2}{8}\left(D_1f(y_1)-D_1f(y_0)\right)  +    \frac{h^3}{48}\left(D_2f(y_1)+D_2f(y_0)\right), $ \\ $\Longleftrightarrow$ 
		
		$ \left\{
		\begin{array}{llll}\displaystyle
			y_{1/2}=y_0+\frac{h}{2}f(y_0) + \frac{h^2}{8} D_1f(y_{0})  + \frac{h^3}{48}D_2f(y_{0}),  \\[.25cm]
			\displaystyle y_{1}=y_{1/2}+\frac{h}{2}f(y_{1}) -  \frac{h^2}{8} D_1f(y_{1})  + \frac{h^3}{48}D_2f(y_{1}). 
		\end{array} \right.
		$
	\end{tabular}
\end{center}

The symbols $D_1$ and $D_2$ denote the first and second order Lie derivatives of the vector field. As shown in \cite{Iavernaro2018_torino}, both integrators are conjugate to a symplectic method up to order six.

In this work, we aim to approximate $D_1$ and $D_2$ using suitable finite-difference schemes. Since these terms are multiplied by $h^2$ and $h^3$, respectively, preserving both the order of accuracy and the symmetry of the method requires employing symmetric finite-difference formulas of at least second-order accuracy. For $r=0,1/2,1$, we adopt the following standard centered differences:
\begin{align}
	\hat{D}_1 f_{n+r} &= \frac{f(y_{n+r+\alpha}) - f(y_{n+r-\alpha})}{2\alpha h}, \label{diffMPa} \\
	\hat{D}_2 f_{n+r} &= \frac{f(y_{n+r+\alpha}) - 2f(y_{n+r}) + f(y_{n+r-\alpha})}{\alpha^2 h^2}, \label{diffMPb}
\end{align}
where $\alpha > 0$ is a parameter, and the auxiliary times are given by $t_{n+r \pm \alpha} = t_{n+r} \pm \alpha h$.

These formulas require the auxiliary local stages $y_{n+r \pm \alpha}$, which must be computed at each step. In \cite{Mathematica2021}, these additional approximations were generated using two distinct approaches: a second-order explicit Runge--Kutta method and the trapezoidal rule. In the latter case, for a suitable value of $\alpha$, a novel symmetric and symplectic fourth-order Runge--Kutta method, along with its conjugate-symplectic twin, was obtained.

This result motivates the question: can an alternative approximation of the auxiliary stages $y_{n+r \pm \alpha}$ lead to the recovery of the Gauss--Legendre integrator of order four and its conjugate-symplectic counterpart? As we shall demonstrate in the following sections, the answer is affirmative when a degree-two collocation polynomial is used.

\section{An MDMP4 variant based on a quadratic collocation polynomial} 
\label{MDMP}

In this section, we approximate the auxiliary stages $y_{n+ 1/2 \pm \alpha}$ by evaluating a collocation polynomial of degree two at the points $t_{n+1/2 \pm \alpha}$. The resulting method, denoting $f(y_r)=f_r$, reads:

\begin{subequations}
	\label{MPCOL}
	\begin{flalign}
		y_{n+\frac{1}{2}-\alpha} &=  y_{n+\frac{1}{2}} - h \alpha \frac{3}{4} \left( f_{n+\frac{1}{2}-\alpha} + \frac{1}{3} f_{n+\frac{1}{2}+\alpha} \right), \label{MPCOL1} \\
		y_{n+\frac{1}{2}+\alpha} &=  y_{n+\frac{1}{2}} + h \alpha \frac{3}{4} \left( \frac{1}{3} f_{n+\frac{1}{2}-\alpha} + f_{n+\frac{1}{2}+\alpha} \right), \label{MPCOL2} \\
		y_{n+\frac{1}{2}} &= \displaystyle y_n + \frac{h}{2} f_{n+1/2}- \frac{h^2}{8} \hat D_1 f_{n+1/2} + \frac{h^3}{48} \hat D_2 f_{n+1/2}, \label{MPCOL3} \\[.25cm]
		y_{n+1} &= \displaystyle y_n + h f_{n+1/2} + \frac{h^3}{24} \hat D_2 f_{n+1/2}. \label{MPCOL4}
	\end{flalign}
\end{subequations}

Expressed as a Runge--Kutta scheme, the method is described by the following tableau:

\begin{equation}\label{MIDPOINT_COLL}
	\begin{array}{c|ccc}
		1/2 - \alpha & -\frac{3}{4} \alpha + \frac{1}{16 \alpha} + \frac{1}{48 \alpha^2} & \frac{1}{2} - \frac{1}{24 \alpha^2} & -\frac{1}{4} \alpha - \frac{1}{16 \alpha} + \frac{1}{48 \alpha^2} \\[.15cm]
		1/2 & \frac{1}{16 \alpha} + \frac{1}{48 \alpha^2} & \frac{1}{2} - \frac{1}{24 \alpha^2} & -\frac{1}{16 \alpha} + \frac{1}{48 \alpha^2} \\[.15cm]
		1/2 + \alpha & \frac{1}{4} \alpha + \frac{1}{16 \alpha} + \frac{1}{48 \alpha^2} & \frac{1}{2} - \frac{1}{24 \alpha^2} & \frac{3}{4} \alpha - \frac{1}{16 \alpha} + \frac{1}{48 \alpha^2} \\[.15cm]
		\hline
		& \frac{1}{24 \alpha^2} & 1 - \frac{1}{12 \alpha^2} & \frac{1}{24 \alpha^2}
	\end{array}
\end{equation}

We denote this method by AMDMP4\_C2 (MDMP4 approximated by a collocation polynomial of degree 2).

\begin{theorem}
	The Runge--Kutta scheme defined by the tableau \eqref{MIDPOINT_COLL} is symplectic if and only if $\alpha = \sqrt{3}/6$. In this case, the method coincides with the Gauss--Legendre method of order 4 and can be obtained as the composition $\Phi_{h/2} \circ \Psi_{h/2}$, where $\Phi_h$ and $\Psi_h$ are defined in \eqref{PhiPsi}.
\end{theorem}

\begin{proof}
	Setting $\alpha = \sqrt{3}/6$, the first and third abscissae coincide with the shifted roots of the Gauss orthogonal polynomial of degree two on the interval $[0,1]$. The coefficient matrix $A$ of the method \eqref{MIDPOINT_COLL} has the second column identically zero, which implies that the nonlinear system associated with the method does not involve the second stage (i.e., the one corresponding to $t_{n+1/2}$). Hence, the second row and column can be removed. A direct computation then shows that the reduced scheme corresponds to the Gauss method of order four, which is known to be symplectic. 
	
	Moreover, $\alpha = \sqrt{3}/6$ is the unique value for which the symplecticity condition 
	$$
	b_i a_{ij} + b_j a_{ji} - b_i b_j = 0, \quad i,j = 1, \dots, 3,
	$$
	is satisfied (see \cite[Theorem 4.3]{HLW06}). Substituting equations \eqref{MPCOL1} and \eqref{MPCOL2} into \eqref{MPCOL3} and \eqref{MPCOL4}, one recovers the two-stage Runge--Kutta methods defined in \eqref{PhiPsi}. \qed
\end{proof}

\smallskip

We emphasize that equations \eqref{MPCOL3} and \eqref{MPCOL4} share the same internal stages, computed via the coupled system \eqref{MPCOL1}--\eqref{MPCOL2}. This means that, when used in composition, only a single nonlinear system must be solved per time step of size $h$. 

\def\bigo{{\cal O}}

We observe that the $1$st and $3$rd internal stages approximate the local solution with an error $\bigo(h^3)$, whereas the
$2$nd stage, $y_{n+1/2}$ approximates it with an error $\bigo(h^4)$.
After evaluating $f(y_{n+1/2})$, we compute $\hat D_1$ and $\hat D_2$ in (\ref{diffMPa})-(\ref{diffMPb})  and we define the polynomial
in the variable $\tau \in [-1/2,1/2]$
\begin{equation}
	\label{dense}
	\displaystyle y_{n+1/2 +  \tau }=y_{n+1/2}+\frac{\tau h}{2}f(y_{n+1/2} ) + \frac{\tau^2 h^2}{2} \hat D_1f(y_{n+1/2})  + \frac{\tau^3 h^3}{6} \hat D_2f(y_{n+1/2}). 
\end{equation}
Note that for $\tau=\pm 1/2$ this dense output passes through the approximations $y_n$ and $y_{n+1}$, so that it provides
a globally continuous approximation of the solution.

\begin{figure}[bht] 
	\begin{center}
		\includegraphics[width=10 cm]{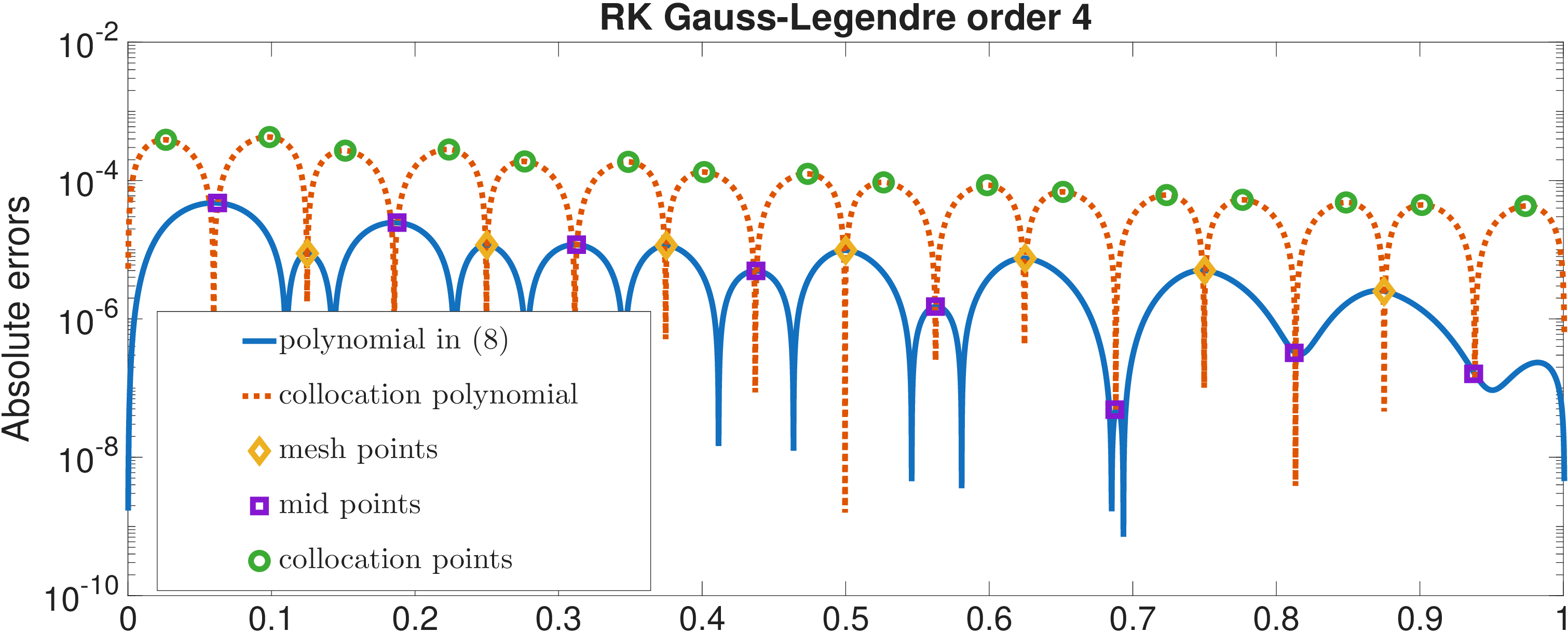}   
	\end{center}
	\caption{Maximum absolute error in the continuous extensions  for the polynomial in  (\ref{dense}) (solid blue line) and the collocation polynomial (dotted red line)).}
	\label{fig:example}
\end{figure}

We recall that the collocation polynomial associated with the Gauss method of order four has degree two and
its error is $\bigo(h^3)$. In contrast, the dense output \eqref{dense} involving discretized derivatives has degree three and its error is
$\bigo(h^4)$. Consequently, with only one additional function evaluation we get a fourth-order continuous extension of the Gauss method.

\begin{table}[hb]  
	\scriptsize
	\begin{center}
		\begin{tabular}{|c|c|c|c|c|}
			\hline
			$h$ & error for the polynomial in (\ref{dense})		& rate & error for the collocation polynomial  & rate \\
			\hline
			1/8&  $4.7402 \cdot 10^{-5} $ &          0  & $4.2624 \cdot 10^{-4}$            & 0\\
			1/16& $3.4239 \cdot 10^{-6} $   &3.79 &  $5.7704 \cdot 10^{-5} $  &2.88\\
			1/32& $2.3042 \cdot 10^{-7} $   &3.89   &$7.4913 \cdot 10^{-6}$   &2.94\\
			1/64& $1.4949 \cdot 10^{-8} $   &3.94 &$9.5368 \cdot 10^{-7} $   &2.97\\
			1/128& $9.5203 \cdot 10^{-10} $   &3.97   &$1.2028 \cdot 10^{-7}$    &2.98\\
			1/256& $6.0064 \cdot 10^{-11} $  &3.98   &$1.5102 \cdot 10^{-8}$    &2.99\\
			1/512 & $3.7718 \cdot 10^{-12}$    &3.99   &$1.8919 \cdot 10^{-9 }$   &2.99\\
			\hline
		\end{tabular}
	\end{center}
	\caption{Maximum absolute error in the continuous extensions for the polynomial in  (\ref{dense}) (continous blue line) and the collocation polynomial (dotted red line).}
	\label{table:example}
\end{table}

	To show the behavior of the derived continuous extension we have solved the following linear boundary value problem:
	\begin{equation}\label{prob1}
		\begin{array}{l} \epsilon y'' = y, \\
			y(0) = 1 , y(1)=0, \end{array}
	\end{equation}
	rewritten as a first order system, whose  exact solution is
	$$\frac{e^{-x/\sqrt{\epsilon}}- e^{ (x-2)/\sqrt{\epsilon}}}{1-e^{-2/\sqrt{\epsilon}}}$$
	The results  using constant stepsize and $\epsilon=0.1$ are reported in  Figure  \ref{fig:example} for $h=1/8$ and Table \ref{table:example}. A comparison with the continuous extension produced by the collocation polynomial associated with the fourth-order Gauss method shows the better performance of formula (\ref{dense}) in terms of order of convergence.

\section{An MDTR4  variant   based on a quadratic collocation polynomial}
\label{MDTR}
We repeat the approximation procedure of the two Lie derivatives appearing in the MDTR4 method: the auxiliary stages are now $y_{n \pm \alpha}$ and  $y_{n+1 \pm \alpha}$. Evaluating a collocation polynomial of degree two at the points $t_{n \pm \alpha}$ and $t_{n+1 \pm \alpha}$ yields the method:
\begin{subequations}
	\label{TRCOL}
	\begin{flalign}
		y_{n+r-\alpha} &=  y_{n+r} - h   \alpha \frac{3}{4}(  f_{n+r-\alpha} +  \frac{1}{3} f_{n+r+\alpha}  ), \label{TRCOL1}\\
		y_{n+r+\alpha} &=  y_{n+r} + h \alpha \frac{3}{4} ( \frac{1}{3}  f_{n+r-\alpha} + f_{n+r+\alpha}) \label{TRCOL2} \\
		y_{n+\frac{1}{2}} &= \displaystyle y_n + \frac{h}{2} f_{n} + \frac{h^2}{8} \hat D_1 f_{n} + \frac{h^3}{48} \hat D_2 f_{n}, \label{TRCOL3} \\
		y_{n+1} &= \displaystyle y_{n+\frac{1}{2}} + \frac{h}{2} f_{n+1} - \frac{h^2}{8} \hat D_1 f_{n+1} + \frac{h^3}{48} \hat D_2 f_{n+1}. \label{TRCOL4}
	\end{flalign}
\end{subequations}

Equations \eqref{TRCOL1} and \eqref{TRCOL2} provide the internal stages of the method \eqref{TRCOL3} for $r=0$ and of the method \eqref{TRCOL4} for $r=1$.  The method resulting from the composition, denoted AMDTR4\_C2 (MDTR4 approximated by a collocation polynomial of degree 2), has the form
\begin{equation}
	\label{AMDTR4C2}
	\begin{array}{rcl}
		y_{n+1}&=&\displaystyle y_n+\frac{h}{2}(f_{n}    +  f_{n+1} )  -  \frac{h^2}{8}\left(\hat D_1f_{n+1}-\hat D_1f_{n} \right)  +  \\[.2cm]
		& & \displaystyle   + \frac{h^3}{48}\left(\hat D_2f_{n+1}+\hat D_2f_{n} \right).
	\end{array}
\end{equation}
This method is of order 4 and is conjugate to the AMDMP4\_C2 method when used with the same parameter $\alpha$. The conjugate symplectic twin associated with the Gauss method of order four is obtained by choosing $\alpha = \sqrt{3}/6$ and is reported in (\ref{PhiPsiConjSymp}). For this value of $\alpha$, the stability function of \eqref{AMDTR4C2} coincides with the Padé rational function of order $(2,2)$. Furthermore, the two methods \eqref{TRCOL3} and \eqref{TRCOL4}, used with step size $h$, correspond to the Runge--Kutta methods $\Psi_h$ and $\Phi_h$ defined in \eqref{PhiPsi}, respectively.

\section{Generalization to higher-order formulae}
\label{sec:gen}

A straight forward extension of the approach of Section~\ref{sec: MDMP_MDTR}, by using higher degree Taylor polynomials,
does not lead to Gauss methods of higher order. We therefore consider an alternative approach.

We start with an $s$-stage Runge--Kutta method with coefficient matrix $A$, distinct quadrature nodes $c = A \one$ and weights $b^\top$.
Here, $\one$ denotes the $s$-dimensional vector with all entries equal to $1$. We then consider the two Runge--Kutta methods:
\begin{description}
	\item[$\Phi_h\,$:]\,%
	with Runge--Kutta matrix $A_1 = 2A$, quadrature nodes $c_1 = 2c$ and weights $b_1^T$, such that the underlying quadrature
	formula is at least
	of order $s$,
	\item[$\Psi_h\,$:]\,%
	with Runge--Kutta matrix $A_2 = 2A - \one b_1^\top $, nodes $c_2 = 2c-\one$ and weights $b_2^\top$ such that the  underlying quadrature
	formula is at least
	of order $s$.
\end{description}
We first note that for the implicit midpoint rule we have $A=(1/2)$, so that $\Phi_h$ becomes the implicit Euler method and $\Psi_h$
the explicit Euler method. For the $2$-stage Gauss method, this gives precisely the two methods of \eqref{PhiPsi}.

	The coefficients of the s-stage Gauss method of order 6 ($s=3$) are tabulated in \cite[page 72]{HW96}.
	The explicit formulas for $\Phi_h$ and $\Psi_h$ then give the desired composition.

\begin{theorem}
	For a given $s$-stage Runge--Kutta method $\Theta_h = (A, c, b^\top )$ of order at least $s$, the composition
	$\Psi_{h/2}\circ \Phi_{h/2}$ is equivalent to $\Theta_h$.
\end{theorem}

\begin{proof}
	The composition $\Psi_{h/2}\circ \Phi_{h/2}$ can be considered as a $2s$-stage Runge--Kutta method with coefficients
	\[
	\Psi_{h/2}\circ \Phi_{h/2} ~~= ~~
	\begin{array}{c|cc}
		c_1/2  & A_1/2 & 0 \\[.15cm]
		\one /2 + c_2/2  &\one b_1^\top /2 & A_2/2 \\[.1cm]
		\hline \\[-3mm]
		& b_1^\top /2 & b_2^\top /2
	\end{array}
	~~= ~~
	\begin{array}{c|cc}
		c  & A & 0 \\[.15cm]
		c  &\one b_1^\top /2 & A -  \one b_1^\top/2 \\[.1cm]
		\hline \\[-3mm]
		& b_1^\top /2 & b_2^\top /2
	\end{array} .
	\]
	This Runge--Kutta method is $S$-reducible (see Definition~IV.12.17 of \cite{HW96}). In fact, the internal stages of the first block
	are identical to the internal stages of the second block. Therefore, this method reduces to
	$(A,c, (b_1+b_2)^\top \! /2 )$. The assumption on the order of the underlying quadrature formulas implies that $(b_1 + b_2)/2 = b$.
	\qed
\end{proof}

\begin{corollary}
	For a symmetric  $s$-stage Runge--Kutta method $\Theta_h = (A, c, b^\top )$ of order at least $s$, $\Phi_h$ and $\Psi_h$ satisfy the relations
	\begin{equation}
		\label{adj}
		\Phi_h^\ast=\Psi_h, \qquad \Psi_h^\ast=\Phi_h,
	\end{equation}
	where $\Phi_h^\ast := \Phi^{-1}_{-h}$ denotes the adjoint method of $\Phi_h$ (and analogously for $\Psi_h^\ast$)  \cite[Definition II.3.1]{HLW06}. 
\end{corollary}
\begin{proof}
	Relations \eqref{adj} can be proven by exploiting the general results provided in \cite[Theorem II.8.3, Theorem II.8.8]{HWN93}. The symmetry of the  method $\Theta_h$ implies that 
	$$ 
	b_j = b_{s+1-j}, \quad a_{i,j}+a_{s+1-i,s+1-j} = b_j, \qquad i,j=1,\dots,s,
	$$  
	and, consequently,  $c_i + c_{s+1-i}=1$.
	We first prove $(b_1)_{i} = (b_2)_{s+1-i}$ for the weights of the methods $\Phi_h$ and $\Psi_h$. These weights are defined by
	$$
	\sum_{i=1}^s (b_1)_i \,(2c_i )^{k-1} = \frac 1k, \qquad \sum_{i=1}^s (b_2)_{s+1-i} \,(1-2c_i)^{k-1} = \frac 1k ,  \qquad k=1,\dots,s,
	$$
	where we have used the relation $2c_{s+1-i} -1 = 1-2c_i$. Using the binomial identity one can prove by induction on $k$ that
	$ \sum_{i=1}^s (b_2)_{s+1-i} \,(2c_i)^{k-1} = \frac 1k $, so that  $(b_1)_{i} = (b_2)_{s+1-i}$.
	
	With this preparation we obtain
	$$
	\begin{array}{rcl}
		(A_2)_{i,j} &=& 2\,a_{i,j}-(b_1)_j = 2\, a_{i,j} - 2\,b_j + (b_2)_j \\[.25cm]
		&=&  (b_2)_j - 2\, a_{s+1-i,s+1-j} = (b_1)_{s+1-j} -  (A_1)_{s+1-i,s+1-j}.
	\end{array}
	$$
	which, together with $(b_2)_j = (b_1)_{s+1-j} $, proves that $\Psi_h$ is the adjoint of $\Phi_h$.
	\qed
\end{proof}

\section{Conclusions} \label{sec:conc}
Starting from a pair of multiderivative Explicit and Implicit Taylor methods, and discretizing the involved derivatives, we have shown that their composition leads to the fourth-order Gauss–Legendre Runge–Kutta method, along with its associated conjugate symplectic scheme.
This construction has naturally provided a higher-order dense output formula, intrinsically linked to the scheme, which requires only a single additional function evaluation.

Analyzing the coefficients of the two composing methods has revealed a way to extend the approach to Gauss–Legendre methods of arbitrary order.

\section*{Acknowledgments}  F. Iavernaro and F. Mazzia  acknowledge the contribution of the National Recovery and Resilience Plan, Mission 4 Component 2 -Investment 1.4 -NATIONAL CENTER FOR HPC, BIG DATA AND QUANTUM COMPUTING -funded by the European Union -NextGenerationEU -(CUP H93C22000450007). F. Iavernaro and F. Mazzia are members of the INdAM GNCS national group and supported by INDAM-GNCS project CUP\_E53C24001950001). The same authors also wish to thank ``{\em Fondo acquisto e manutenzione attrezzature per la ricerca}'' D.R. 3191/2022, Universit\`a degli Studi di Bari Aldo Moro.

\bibliographystyle{elsarticle-num-names}
\bibliography{AMDMPTR}

\end{document}